\documentclass{amsart}
\usepackage{amssymb}

\newtheorem{Theorem}{Theorem}

\def\Z{\mathbb Z}
\def\Q{\mathbb Q}
\def\R{\mathbb R}
\def\H{\mathbb H}
\def\C{\mathbb C}

\def\Im{\text{\rm Im\,}}

\begin{document}

\title{Defining equations of $X_0(2^{2n})$}
\author{Fang-Ting Tu}
\address{Department of Applied Mathematics \\
  National Chiao Tung University \\
  Hsinchu 300 \\
  TAIWAN}
\email{ft.am95g@cc.nctu.edu.tw}
\author{Yifan Yang}
\address{Department of Applied Mathematics \\
  National Chiao Tung University \\
  Hsinchu 300 \\
  TAIWAN}
\email{yfyang@math.nctu.edu.tw}
\date{31 January 2007}
\thanks{The authors were support by Grant 95-2115-M-009-005 of the
  National Science Council (NSC) of Taiwan.}

\begin{abstract} In this note we will obtain defining equations of
  modular curves $X_0(2^{2n})$. The key ingredient is a recursive
  formula for certain generators of the function fields on $X_0(2^{2n})$.
\end{abstract}
\maketitle

\begin{section}{Introduction and statements of results}
Let $\Gamma$ be a congruence subgroup of $SL_2(\R)$ commensurable with
$SL_2(\Z)$. The modular curve $X(\Gamma)$ is defined as the quotient
of the extended upper half-plane
$\H^\ast=\{\tau\in\C:\Im\tau>0\}\cup\mathbb P^1(\Q)$ by the action
of $\Gamma$. It has a complex structure as a compact Riemann surface
(i.e., a non-singular irreducible projective algebraic curve),
and the polynomials defining the Riemann surface are called
{\it defining equations} of $X(\Gamma)$. The problem of explicitly
determining the equations of modular curves has been addressed by
numerous authors. For instance, Galbraith \cite{Galbraith1},
Murabayashi \cite{Murabayashi}, and Shimura \cite{ShimuraM} used the
so-called canonical embeddings to find equations of $X_0(N)$ that are
non-hyperelliptic. For hyperelliptic $X_0(N)$, we have results of
Galbraith \cite{Galbraith1}, Gonz\'alez \cite{Gonzalez}, Hibino
\cite{Hibino}, Hibino-Murabayashi \cite{Hibino-Murabayashi}, and
Shimura \cite{ShimuraM}. In \cite{Reichert} Reichert used the fact
that $X_1(N)=X(\Gamma_1(N))$ is the moduli space of isomorphism
classes of elliptic curves with level $N$ structure to compute
equations of $X_1(N)$ for $N=11,~13,\ldots,18$. Furthermore, in
\cite{Ishida-Ishii} Ishida and Ishii proved that for each $N$ two
certain products of the Weierstrass $\sigma$-functions generate the
function field on $X_1(N)$, and thus the relation between these two
functions defines $X_1(N)$. A similar method was employed in
\cite{Ishida} to obtain equations of $X(N)=X(\Gamma(N))$. Very
recently, in \cite{Yang2} the second author of the present article
devised a new method for obtaining defining equations of $X_0(N)$,
$X_1(N)$, and $X(N)$, in which the required modular functions are
constructed using the generalized Dedekind eta functions. (See
\cite{Yang} for the definition and properties of these functions.)

When $\Gamma_1$ and $\Gamma_2$ are two congruence subgroups such that
$\Gamma_2$ is contained in $\Gamma_1$ and a defining equation of
$X(\Gamma_1)$ is known, one may attempt to deduce an equation for
$X(\Gamma_2)$ using the natural covering $X(\Gamma_2)\to
X(\Gamma_1)$. Of course, the main difficulty in this approach lies at
finding an explicit description of the covering map. In this note we
will prove a recursive formula for the coverings
$X_0(2^{2(n+1)})\to X_0(2^{2n})$, from which we easily obtain
defining equations of $X_0(2^{2n})$ for positive integers $n$.

To state our result, we first recall the definition of the Jacobi
theta functions
$$
  \theta_2(\tau)=\sum_{n\in\Z}q^{(2n+1)^2/8}
 =2\frac{\eta(2\tau)^2}{\eta(\tau)},
$$
$$
  \theta_3(\tau)=\sum_{n\in\Z}q^{n^2/2}
 =\frac{\eta(\tau)^5}{\eta(\tau/2)^2\eta(2\tau)^2},
$$
and
$$
  \theta_4(\tau)=\sum_{n\in\Z}(-1)^nq^{n^2/2}
 =\frac{\eta(\tau)^2}{\eta(\tau/2)^2},
$$
where $q=e^{2\pi i\tau}$ and
$$
  \eta(\tau)=q^{1/24}\prod_{n=1}^\infty(1-q^n)
$$
is the Dedekind eta function. Now our main result can be stated as
follows.

\begin{Theorem} Let $P_6(x,y)=y^4-x^3-4x$, and for $n\ge 7$ define
  polynomials $P_n(x,y)$ recursively by
$$
  P_n(x,y)
 =\left[Q_{n-1}\left(\frac{\sqrt{x^2+4}}{\sqrt x},\frac y{\sqrt x}\right)^2
 -R_{n-1}\left(\frac{\sqrt{x^2+4}}{\sqrt x},\frac y{\sqrt
  x}\right)^2\right] x^{2^{n-5}},
$$
where $Q_{n-1}(x,y)=(P_{n-1}(x,y)+P_{n-1}(-x,y))/2$ and
$R_{n-1}=P_{n-1}-Q_{n-1}$. Then $P_{2n}(x,y)=0$ is a defining equation
of the modular curve $X_0(2^{2n})$ for $n\ge 3$.

To be more precise, for $n\ge 6$, let
$$
  x_n=\frac{2\theta_3(2^{n-1}\tau)}{\theta_2(2^{n-1}\tau)}, \qquad
  y_n=\frac{\theta_2(8\tau)}{\theta_2(2^{n-1}\tau)}.
$$
Then
\begin{enumerate}
\item $x_{n-1}=\sqrt{(x_n^2+4)/x_n}$ and
  $y_{n-1}=y_n/\sqrt{x_n}$.
\item $P_{n}(x_n,y_n)=0$, and $P_n(x,y)$ is irreducible over $\Q$.
\item When $n$ is an even integer greater than $4$, $x_n$ and $y_n$ are
  modular functions on $\Gamma_0(2^n)$ that are holomorphic
  everywhere except for a pole of order $2^{n-4}$ and $2^{n-4}-1$,
  respectively, at $\infty$. (Thus, they generate the field of modular
  functions on $\Gamma_0(2^n)$ and the relation $P_n(x_n,y_n)=0$
  between them is a defining equation for $X_0(2^n)$.) 
\end{enumerate}
\end{Theorem}

We remark that from the definition of $Q_n$ and $R_n$, it is easy to
see that $Q_{n-1}(x,y)^2-R_{n-1}(x,y)^2$ is a polynomial of degree
$2^{n-4}$ contained in $\Z[x^2,y^4]$. Thus, $P_n(x,y)$ are indeed
polynomials. We also remark that when $n$ is odd, the polynomial
$P_n(x,y)$ fails to be a defining equation of $X_0(2^n)$ because in
this case
$$
  y_n(\tau)=\frac{\eta(16\tau)^2\eta(2^{n-1}\tau)}
  {\eta(8\tau)\eta(2^n\tau)^2}
$$
is not modular on $\Gamma_0(2^n)$. (When $n$ is odd, $y_n$ does not
satisfy the conditions of Newman \cite[Theorem I]{Newman} for a
product of Dedekind eta functions to be modular on $\Gamma_0(N)$.
Indeed, one can show that when $n$ is odd,
$$
  y_n\left(\frac{a\tau+b}{c\tau+d}\right)
 =\left(\frac 2d\right)y_n(\tau), \qquad
  \begin{pmatrix}a&b\\ c&d\end{pmatrix}\in \Gamma_0(2^n),
$$
where $\left(\frac\cdot d\right)$ is the Jacobi symbol.)
\medskip

\noindent{\bf Examples.} Using Theorem 1, we find that a defining
equation of $X_0(256)$ is
$$
  y^{16}-16x(x+2)^4(x^2+4)y^8-x(x+2)^4(x-2)^8(x^2+4)=0,
$$
and an equation for $X_0(1024)$ is
\begin{equation*}
\begin{split}
 &y^{64}-4096vy^{56}-61696uvy^{48}-512uv(253u+30464v)y^{40} \\
 &\qquad-16uv(4619u^2-2^8\cdot2053uv
  +2^{16}\cdot7\cdot73v^2)y^{32} \\
 &\qquad-512uv(31u^3+35712u^2v+3\cdot2^{16}uv^2+2^{23}v^3)y^{24} \\
 &\qquad-32u^3v(47u^2-320000uv+2^{15}\cdot17\cdot31v^2)y^{16} \\
 &\qquad-64u^3v(u^3+5248u^2v+2^{18}\cdot 5uv^2+2^{26}v^3)y^8-u^7v=0,
\end{split}
\end{equation*}
where $u=(x-2)^8$ and  $v=x(x+2)^4(x^2+4)$.
\medskip

Our interest in the modular curves $X_0(2^{2n})$ stems from the
following remarkable observation of Hashimoto. When $n=3$,
it is known that the curve $X_0(64)$ is non-hyperelliptic (see
\cite{Ogg-hyperelliptic}) of genus $3$. Then the theory of Riemann
surfaces says that it can be realized as a plane quartic. Indeed, it
can be shown that the space of cusp forms of weight $2$ on
$\Gamma_0(64)$ is spanned by
$$
  x=\eta(4\tau)^2\eta(8\tau)^2, \qquad
  y=2\eta(8\tau)^2\eta(16\tau)^2, \qquad
  z=\frac{\eta(8\tau)^8}{\eta(4\tau)^2\eta(16\tau)^2},
$$
and the map $X_0(64)\to\mathbb P^2(\C)$ defined by
$\tau\mapsto[x(\tau):y(\tau):z(\tau)]$ is an embedding. Then the
relation
$$
  x^4+y^4=z^4
$$
among $x$, $y$, $z$ is a defining equation of $X_0(64)$ in $\mathbb
P^2$. (The Fermat curve $X^4+Y^4=1$ is birationally equivalent to
$y^4-x^3-4x=0$ in Theorem 1 via the map
$$
  X=\frac{x-2}{x+2}, \qquad Y=\frac{2y}{x+2}.
$$)
Then Hashimoto pointed out the curious fact that the Fermat curve
$F_{2^n}:x^{2^n}+y^{2^n}=1$ and the modular curve $X_0(2^{2n+2})$ have
the same genus for all positive integer $n$. In fact, there are more
similarities between these two families of curves. For instance, the
obvious covering $F_{2^{n+1}}\to F_{2^n}$ given by
$[x:y:z]\to[x^2:y^2:z^2]$ branches at $3\cdot 2^n$ points, each of
which is of order 2. On the other hand, the congruence
subgroup $\Gamma_0(2^{2n+2})$ is conjugate to
$$
  \Gamma_0^0(2^{n+1})=\left\{\begin{pmatrix}a&b\\ c&d\end{pmatrix}
  \in SL_2(\Z):2^{n+1}|b,c\right\},
$$
and the natural covering $X_0^0(2^{n+2})\to X_0^0(2^{n+1})$ also branches
at $3\cdot 2^n$ cusps of $X_0^0(2^{n+1})$. These observations
naturally lead us to consider the problem whether the modular curve
$X_0(2^{2n+2})$ is birationally equivalent the Fermat curve $F_{2^n}$.
It turns out that this problem can be answered easily as follows.

According to \cite{Elkies, Kenku-Momose, Ogg-automorphisms}, when a
modular curve $X_0(N)$ has genus $\ge 2$, any automorphism of $X_0(N)$
will arise from the normalizer of $\Gamma_0(N)$ in $SL_2(\R)$, with
$N=37,~63$ being the only exceptions. Now by
\cite[Theorem 8]{Atkin-Lehner}, for all $n\ge 7$, the index of
$\Gamma_0(2^n)$ in its normalizer in $SL_2(\R)$ is $128$.
% (When $n\ge
%7$, the index $128$ is accounted for by the coset representatives
%$\begin{pmatrix}1&1/8\\ 0&1\end{pmatrix}^i
% \begin{pmatrix}1&0\\ 2^{n-3}&1\end{pmatrix}^jW_2^k$, where
%$i,j=0,\ldots,7$, $k=0,1$, and $W_2$ is the Atkin-Lehner involution.)
Therefore, the automorphism group of $X_0(2^{2n+2})$ has order $128$
for all $n\ge 3$. On the other hand, it is clear that the automorphism
group of any Fermat curve contains $S_3$. Thus, we conclude that the
modular curve $X_0(2^{2n+2})$ cannot be birationally equivalent to the
Fermat curve $F_{2^n}$ when $n\ge 3$. Still, it would be an
interesting problem to study the exact relation between these two
families of curves.
\medskip

\centerline{\sc Acknowledgment}
\medskip

The authors would like to thank Professor Hashimoto of the Waseda
University for drawing their attention to the family of modular curves
$X_0(2^n)$ and for several enlightening conversations. The authors would
also like to thank Professor M. L. Lang of the National University of
Singapore for providing information about normalizers of congruence
subgroups.
\end{section}

\begin{section}{Proof of Theorem 1} To prove
$x_{n-1}=\sqrt{(x_n^2+4)/x_n}$, we first verify the case $n=2$ by
comparing the Fourier expansions for enough terms, and then the
general case follows since $x_n(\tau)$ is actually equal to
$x_1(2^{n-1}\tau)$. The proof of $y_{n-1}=y_n/\sqrt{x_n}$ is equally
simple. We have
$$
  \frac{y_{n-1}^2}{y_n^2}=\frac{\theta_2(2^{n-1}\tau)^2}
 {\theta_2(2^{n-2}\tau)^2}=\frac
 {\eta(2^{n-2}\tau)^2\eta(2^n\tau)^4}{\eta(2^{n-1}\tau)^6}
 =\frac{\theta_2(2^{n-1}\tau)}{\theta_3(2^{n-1}\tau)}
 =\frac1{x_n}.
$$
This proves the recursion part of the theorem. We now show that when
$n\ge 6$ is an even integer, $x_n$ and $y_n$ are modular functions on
$\Gamma_0(2^n)$ that have a pole of order $2^{n-4}$ and $2^{n-4}-1$,
respectively, at $\infty$ and are holomorphic everywhere.

By the criteria of Newman \cite{Newman}, a product
$$
  \prod_{k=0}^n\eta(2^k\tau)^{e_k}
$$
of Dedekind eta functions is a modular function on $\Gamma_0(2^n)$ if
the four conditions
\begin{enumerate}
\item $\sum_k e_k=0$,
\item $\sum_k ke_k\equiv 0\text{ mod }2$,
\item $\sum_k e_k2^k\equiv 0\text{ mod }24$,
\item $\sum_k e_k2^{n-k}\equiv 0\text{ mod }24$,
\end{enumerate}
are satisfied. Now we have
$$
  x_n=\frac{\eta(2^{n-1}\tau)^6}{\eta(2^{n-2}\tau)^2\eta(2^n\tau)^4},
  \qquad
  y_n=\frac{\eta(16\tau)^2\eta(2^{n-1}\tau)}
  {\eta(8\tau)\eta(2^n\tau)^2}.
$$
It is clear that when $n$ is an even integer greater than $2$, the
four conditions are all satisfied for $x_n$ and $y_n$. We now show
that $x_n$ and $y_n$ have poles only at $\infty$ of the claimed order.

Still assume that $n\ge 4$ is an even integer. Since $x_n$ and $y_n$
are $\eta$-products, they have no poles nor zeros in $\H$. Also, it
can be checked directly that $x_n$ and $y_n$ have a pole of order
$2^{n-4}$ and $2^{n-4}-1$, respectively, at $\infty$. It remains to
consider other cusps. For an odd integer $a$ and
$k\in\{0,1,\ldots,n-1\}$, the width of the cusp $a/2^k$ is
$$
  h_{n,k}=\begin{cases}1, &\text{if }k\ge n/2, \\
  2^{n-2k}, &\text{if }k<n/2.
  \end{cases}
$$
Choosing a matrix $\sigma=\begin{pmatrix}a&b\\ 2^k&d\end{pmatrix}$ in
$SL_2(\Z)$, a local parameter at $a/2^k$ is
$$
  e^{2\pi i\sigma^{-1}\tau/h_{n,k}}.
$$
Therefore, the order of a function $f(\tau)$ at $a/2^k$ is the same as
the order of $f(\sigma\tau)$ at $\infty$, multiplied by $h_{n,k}$.

Now recall that, for $\alpha=\begin{pmatrix}a&b\\ c&d\end{pmatrix}\in
SL_2(\Z)$, we have
$$
  \theta_2(\tau)\big|\alpha=
  \begin{cases}
  \epsilon q^{1/8}+\cdots, &\text{if }2|c, \\
  \epsilon+\cdots, &\text{if }2\nmid c,
  \end{cases}
$$
and
$$
  \theta_3(\tau)\big|\alpha=
  \begin{cases} \epsilon+\cdots, &\text{if }2|ac, \\
  \epsilon q^{1/8}+\cdots, &\text{if }2\nmid ac,
  \end{cases}
$$
where $\epsilon$ represents a complex number, but may not be the same at
each occurence. (Up to multipliers, if $\alpha$ is congruent to the
identity matrix or $\begin{pmatrix}1&1\\ 0&1\end{pmatrix}$ modulo $2$,
then the action of $\alpha$ fixes $\theta_2$. Any other matrices will
send $\theta_2$ to either $\theta_3$ or $\theta_4$. This explains the
fact about $\theta_2$. The fact about $\theta_3$ can be explained
similarly.) When $k=n-1$, we have
$$
  2^{n-1}\begin{pmatrix}a&b\\ 2^{n-1}&d\end{pmatrix}\tau
 =\frac{a(2^{n-1}\tau)+2^{n-1}b}{(2^{n-1}\tau)+d}
 =\begin{pmatrix}a& 2^{n-1}b\\ 1&d\end{pmatrix}(2^{n-1}\tau)
$$
and
$$
  8\begin{pmatrix}a&b\\ 2^{n-1}&d\end{pmatrix}\tau
 =\frac{a(8\tau)+8b}{2^{n-4}(8\tau)+d}
 =\begin{pmatrix}a& 8b\\ 2^{n-4}&d\end{pmatrix}(8\tau).
$$
It follows that
$$
  x_n\left(\begin{pmatrix}a&b\\ 2^{n-1}&d\end{pmatrix}\tau\right)
 =\frac{\epsilon_1q^{2^{n-4}}+\cdots}{\epsilon_2+\cdots}
 =\epsilon q^{2^{n-4}}+\cdots,
$$
and
$$
  y_n\left(\begin{pmatrix}a&b\\ 2^{n-1}&d\end{pmatrix}\tau\right)
 =\frac{\epsilon_1 q+\cdots}{\epsilon_2+\cdots}
  =\epsilon q+\cdots.
$$
That is, $x_n$ and $y_n$ have a zero of order $2^{n-4}$ and $1$,
respectively, at $a/2^{n-1}$.

When $k=4,\ldots,n-2$, we have
$$
  2^{n-1}\begin{pmatrix}a&b\\ 2^k&d\end{pmatrix}\tau
 =\begin{pmatrix}2^{n-k-1}a&-1\\ 1&0\end{pmatrix}
  (2^{2k-n+1}\tau+2^{k-n+1}d),
$$
$$
  8\begin{pmatrix}a&b\\ 2^k&d\end{pmatrix}\tau
 =\frac{a(8\tau)+8b}{2^{k-3}(8\tau)+d}
 =\begin{pmatrix}a& 8b\\ 2^{k-3}&d\end{pmatrix}(8\tau).
$$
Therefore,
$$
  x_n\left(\begin{pmatrix}a&b\\ 2^k&d\end{pmatrix}\tau\right)
 =\frac{\epsilon_1+\cdots}{\epsilon_2+\cdots}=\epsilon+\cdots,
$$
and
$$
  y_n\left(\begin{pmatrix}a&b\\ 2^k&d\end{pmatrix}\tau\right)
 =\frac{\epsilon_1 q+\cdots}{\epsilon_2+\cdots}
  =\epsilon q+\cdots.
$$
In other words, $x_n$ has no poles nor zeros at $a/2^k$ for
$k=4,\ldots,n-4$, while $y_n$ has zeros of order $h_{n,k}$ at those
points.

When $k=0,\ldots,3$, we have
$$
  2^{n-1}\begin{pmatrix}a&b\\ 2^k&d\end{pmatrix}\tau
 =\begin{pmatrix}2^{n-k-1}a&-1\\ 1&0\end{pmatrix}
  (2^{2k-n+1}\tau+2^{k-n+1}d),
$$
$$
  8\begin{pmatrix}a&b\\ 2^k&d\end{pmatrix}\tau
 =\begin{pmatrix}2^{3-k}a&-1\\ 1&0\end{pmatrix}
  (2^{2k-3}\tau+2^{k-3}d),
$$
and we find that $x_n$ and $y_n$ have no zeros nor poles at $a/2^k$,
$k=0,\ldots,3$.

In summary, we have shown that $x_n$ and $y_n$ have a
pole of order $2^{n-4}$ and $2^{n-4}-1$, respectively, at $\infty$ and
are holomorphic at any other points. Since $2^{n-4}$ and $2^{n-4}-1$
are clearly relatively prime, $x_n$ and $y_n$ generate the field of
modular functions on $X_0(2^n)$. It remains to show that
$P_n$ is irreducible over $\Q$ and $P_n(x_n,y_n)=0$.

When $n=6$, we verify by a direct computation that
$y_6^4-x_6^3-4x_6=0$. Then the recursive formulas for $x_n$ and $y_n$
implies that $P_n(x_n,y_n)=0$ for all $n\ge 6$. Finally, by the theory
of algebraic curve (see \cite[p.194]{Fulton}), the field of modular
functions on $X_0(2^n)$ is an extension field of $\C(x_n)$ of degree
$2^{n-4}$. In other words, the minimal polynomial of $y_n$ over
$\C(x_n)$ has degree $2^{n-4}$. Now it is easy to see that
$P_n(x,y)=0$ is a polynomial of degree $2^{n-4}$ in $y$ with leading
coefficient $1$. We therefore conclude that $P_n$ is irreducible. This
completes the proof of Theorem 1.
\end{section}

\bibliographystyle{plain}

\end{document}